\documentclass[11pt]{amsart}
\usepackage{amssymb}

\begin{document}

\parskip=\smallskipamount

\newtheorem{theorem}{Theorem}[section]
\newtheorem{lemma}[theorem]{Lemma}
\newtheorem{corollary}[theorem]{Corollary}
\newtheorem{proposition}[theorem]{Proposition}

\theoremstyle{definition}
\newtheorem{definition}[theorem]{Definition}
\newtheorem{example}[theorem]{Example}
\newtheorem{xca}[theorem]{Exercise}
\newtheorem{problem}[theorem]{Problem}
\newtheorem{remark}[theorem]{Remark}
\newtheorem{question}[theorem]{Question}
\newtheorem{conjecture}[theorem]{Conjecture}

\newcommand{\B}{\mathbb{B}}
\newcommand{\C}{\mathbb{C}}
\newcommand{\D}{\mathbb{D}}
\newcommand{\N}{\mathbb{N}}
\newcommand{\Q}{\mathbb{Q}}
\newcommand{\Z}{\mathbb{Z}}
\renewcommand{\P}{\mathbb{P}}
\newcommand{\R}{\mathbb{R}}
\newcommand{\T}{\mathbb{T}}
\newcommand{\U}{\mathbb{U}}

\newcommand{\cA}{\mathcal{A}}
\newcommand{\cB}{\mathcal{B}}
\newcommand{\cR}{\mathcal{R}}
\newcommand{\cC}{\mathcal{C}}
\newcommand{\cD}{\mathcal{D}}
\newcommand{\cF}{\mathcal{F}}
\newcommand{\cH}{\mathcal{H}}
\newcommand{\cL}{\mathcal{L}}
\newcommand{\cN}{\mathcal{N}}
\newcommand{\cO}{\mathcal{O}}
\newcommand{\cS}{\mathcal{S}}
\newcommand{\cZ}{\mathcal{Z}}
\newcommand{\cP}{\mathcal{P}}
\newcommand{\cT}{\mathcal{T}}
\newcommand{\cU}{\mathcal{U}}
\newcommand{\cX}{\mathcal{X}}
\newcommand{\cW}{\mathcal{W}}

\newcommand\wt{\widetilde}
\newcommand\spsh{strongly plurisubharmonic}
\newcommand\hra{\hookrightarrow}
\newcommand\End{\mathrm{End}}
\newcommand\Hom{\mathrm{Hom}}
\newcommand\Mono{\mathrm{Mono}}

\def\di{\partial}
\def\dibar{\overline\partial}
\def\bs{\backslash}
\def\wh{\widehat}

\def\a{\alpha}
\def\e{\epsilon}
\def\b{\beta}
\def\l{\lambda}

\def\rank{\mathrm{rank}\,}

\numberwithin{equation}{section}
\setcounter{section}{0}

\title[Approximation of holomorphic maps] 
{Approximation of holomorphic maps with a lower bound on the rank}
\author{Dejan Kolari\v{c}}
\address{Institute of Mathematics, Physics and Mechanics, 
University of Ljubljana, Jadranska 19, 1000 Ljubljana, Slovenia}
\email{dejan.kolaric@fmf.uni-lj.si}

\thanks{Work on this paper was supported by ARRS, Republic of Slovenia.}

%
%    General info
%
\subjclass[2000]{32E30, 32H02, 32M17, 32Q28}
\date{\today} 
\keywords{Holomorphic maps, approximation, transversality, algebraic sets}

\begin{abstract}
Let $K$ be a closed polydisc or ball in $\C^n$, and let $Y$ be a quasi 
projective algebraic manifold which is Zariski locally equivalent to $\C^p$, 
or a complement of an algebraic subvariety of codimension $\ge 2$ 
in such manifold. If $r$ is an integer satisfying $(n-r+1) (p-r+1)\geq 2$
then every holomorphic map from a neighborhood of $K$ to $Y$
with rank $\ge r$ at every point of $K$ can be approximated 
uniformly on $K$ by entire maps $\C^n\to Y$ with rank $\ge r$
at every point of $\C^n$.
\end{abstract}
\maketitle

\section{Introduction}

In this paper we consider the following problem of approximating holomorphic maps
with a lower bound on their rank. Let $K$ be a closed polydisc (or a closed ball) 
in a complex Euclidean space $\C^n$, and let $f$ be a holomorphic map 
from an open neighborhood of $K$ to a complex manifold $Y$ 
such that ${\rm rank}_z f :=\rank (df_z) \ge r$ for all $z\in K$, where $r$ 
is an integer satisfying $1\le r\le \min\{n,\dim Y\}$.
{\em Is it possible to approximate $f$ uniformly on $K$ by entire maps 
$\wt f\colon \C^n\to Y$ satisfying $\rank_z \wt f\ge r$ at every point 
$z\in \C^n$?} 

The answer clearly depends on the complex analytic properties of $Y$.
If $Y$ is Kobayashi hyperbolic \cite{Ko1}, this fails already
when $r=1$ and $K$ is a disc in $\C$. More generally, if $Y$ is Eisenman
$k$-hyperbolic for some $1\le k\le \dim Y$ \cite{Ei} then $Y$ admits no
holomorphic maps $\C^n\to Y$ of rank $\ge k$, 
and hence the answer is negative for $r\ge k$. 
More precise quantitative obstructions to  the existence of large 
polydiscs in complex manifolds were obtained by Kodaira \cite{Kod}.

In the positive direction,  Forster proved  that holomorphic maps 
$\C^n\to Y=\C^p$ satisfy the jet transversality theorem \cite{Fo}, 
which gives a positive answer if $r$ is sufficiently small
compared to $n$ and $p$ (see Theorem \ref{generic1} below). 
If  $r=n<p$, the above rank condition is satisfied by immersions $\C^n\to\C^p$, 
and in this case an affirmative answer follows from the $h$-principle 
due to Eliashberg and Gromov \cite{GR1}.   
If $n>r=p$, the rank condition is satisfied by submersions,
and the approximation result follows from the $h$-principle 
proved by Forstneri\v c \cite{F3}. The Runge approximation problem 
for holomorphic immersions $\C^n\to\C^n$ in the equidimensional case is still open.

In this paper we consider maps to certain algebraic manifolds.
We shall say that a $p$-dimensional complex manifold is of {\em Class $\cA_0$} 
if it is quasi projective algebraic and is covered by finitely 
many Zariski open sets biregularly isomorphic to $\C^p$.
Examples include all complex projective spaces and Grassmanians.
{\em Class $\cA$} will consist of all algebraic manifolds of the form 
$Y=\wh Y\bs A$ where $\wh Y$ is a manifold of class $\cA_0$ and
$A$ is a closed algebraic subvariety of $\wh Y$ of complex codimension
at least two. (See Def.\ \ref{classA} in \S 2.)
The following is our main result.

\begin{theorem}
\label{theorem1}
Let $K\subset \C^n$ be a closed polydisc, a closed ball, or a product of a 
(lower dimensional) closed polydisc and a ball.
Let $Y$ be a $p$-dimensional manifold of Class $\cA$. Assume that $1\le r\le \min\{n,p\}$,
and that $r<n$ if  $n=p$. Every holomorphic map $f$ from a neighborhood of 
$K$ to $Y$ and satisfying $\rank_z f \ge r$ at every point $z\in K$ can be approximated 
uniformly on $K$ by entire maps $\C^n\to Y$ with rank $\ge r$
at every point of $\C^n$.
\end{theorem}

We wish to emphasize that Theorem \ref{theorem1} does not follow 
from the jet transversality theorem, except for values of $r$ which are
small compared to $n$ and $p$; compare with Theorem \ref{generic1} below. 
The following special case may be of particular interest; the analogous
result for submersions (when $n>\dim Y$) was proved in \cite{F3}.

\begin{corollary}
\label{cor1}
Let $K\subset\C^n$ and $Y$ be as in Theorem \ref{theorem1}.
If $n<\dim Y$ then every holomorphic immersion from a neighborhood of 
$K$ to $Y$  can be approximated uniformly on $K$ by 
entire immersions $\C^n\to Y$.
\end{corollary}

Manifolds of Class ${\cA}$ were considered by Gromov under the name 
Ell-regular manifolds \cite[\textsection 3.5]{GR}, and 
by Forstneri\v{c} \cite[\textsection 2]{AJM}; 
in those papers the reader can find many further examples. 
Both classes are stable with respect to blowing up points.
Every such manifold $Y$ enjoy the following properties 
which will play an important role in our proof:
\begin{itemize}
\item[---] holomorphic maps from (neighborhoods of) compact convex sets in $\C^n$
to $Y$ can be approximated by regular algebraic maps (morphisms) $\C^n \to Y$ 
\cite[Corollary 1.2]{AJM};
\item[---] algebraic maps from affine algebraic manifolds to $Y$ 
enjoy a version of the jet transversality theorem (see \cite[Sect.\ 5]{AJM} and
\S 3 below). 
\end{itemize}

\begin{example}
\label{hyperbol1}
Theorem \ref{theorem1} fails in general for maps to manifolds of the form
$Y=\wh Y\backslash A$ where $\wh Y$ is of Class $\cA_0$ and $A$ is a 
complex analytic (but not algebraic) subvariety of codimension at least two,
or if $A$ contains a hypersurface component. We recall a few examples to this effect:

(1) In \cite[\S 5]{RR} Rosay and Rudin constructed discrete sets $D\subset \C^n$
such that the only holomorphic map $F \colon \C^n\to \C^n$ with non degenerate 
Jacobian ($J F\not \equiv 0$) and satisfying $F(\C^n\backslash D)\subset \C^n\backslash D$ 
is the identity map $F(z)=z$ $(z\in \C^n)$; thus any holomorphic map  
$F\colon \C^n\to \C^n\backslash D$ has ${\rm rank}\; F <n$ at each point. 

(2) In \cite[\S 6]{F1} a proper holomorphic embedding $F\colon \C^m\hookrightarrow \C^{m+n}$
is constructed for any pair of integers $m,n\in\N$ such that the image of 
any holomorphic map $G \colon \C^n\to \C^{m+n}$ 
satisfying ${\rm rank}_z G=n$ at some point $z\in\C^n$ intersects 
the submanifold $A=F(\C^m) \subset \C^{m+n}$ infinitely many times.
It follows that any entire map 
$\C^n\to \C^{m+n}\backslash A$ has ${\rm rank}<n$ at each point. 

(3) According to Corollary 2 in \cite{DEG} the complement $\P_2\backslash A$ of a 
very generic algebraic curve $A$ of degree $d\geq 21$ in the projective plane $\P_2$ is 
hyperbolic, i.e., there exist no nonconstant holomorphic maps 
$\C\to \P_2\backslash A$, and hence Theorem \ref{theorem1} fails for $r=1$.
\end{example}

We now give another result whose 
main ingredient is the jet transversality theorem for holomorphic maps.
In many analytic applications it is important to know the dimension 
of `degeneration sets' of a generically chosen holomorphic map $f\colon X\to Y$ 
between a given pair of complex manifolds. Denote by $\cH(X,Y)$ the space of holomorphic 
maps $X\to Y$ equiped with the compact-open topology. By $J^k(X,Y)$ we denote 
the manifold of all $k$-jets of holomorphic maps $X\to Y$. 
Given $f\in \cH(X,Y)$ and an integer $r\in \N$ we set
\[
	\Sigma_{f,r}=\{x\in X \colon {\rm rank}_x f <r \}.
\]

We shall say that $\rank f\ge r$ on a set $K\subset X$ if
$\rank_x f \ge r$ for all $x\in K$; if we do not specify $K$,
it will be understood that $K=X$.

Following \cite{ANN} we say that a complex manifold $Y$ satisfies 
the {\em Convex Approximation Property} (CAP) if every holomorphic map 
from a neighborhood of a compact convex set $K\subset \C^m$ $(m\in \N)$ 
to $Y$ can be approximated uniformly on $K$ by entire maps $\C^m\to Y$. 
By the main result of \cite{ANN} CAP is equivalent to the classical 
Oka property. Examples of complex manifolds with CAP include complex Lie groups
and complex homogeneous spaces.

\begin{theorem}
\label{generic1}
Let $X$ be a Stein manifold and let $Y$ be a complex manifold satisfying 
{\rm CAP}.  Let $\dim X=n$, $\dim Y = p$, and let $r$ be an integer
satisfying $r\leq \min (n,p)$. Set $d=(n-r+1)(p-r+1)$.
\begin{itemize}
\item[(1)] If $d\le n$ then the set 
$\Omega=\{f\in \cH(X,Y) \colon \dim \Sigma_{f,r} = n-d\}$ is open and 
everywhere dense in $\cH(X,Y)$.
\item[(2)] If $d>n$ then the set $\Omega'=\{f\in \cH(X,Y) \colon {\rm rank} f \geq r\}$ 
is open and everywhere dense in $\cH(X,Y)$.
\end{itemize}
\end{theorem}

\begin{remark}
The examples in \ref{hyperbol1} show that Theorem \ref{generic1} fails in general 
if we do not assume anything on $Y$. In the proof we shall use the
jet transversality theorem for holomorphic maps $X\to Y$ which holds 
if $X$ is a Stein manifold and $Y$ satisfies CAP
\cite[Theorem 1.4]{AJM} or the Ell$_1$ property introduced by 
Gromov (see Definition \ref{Ell1} below and \cite[Theorem 4.2]{AJM}). 

On the other hand, Kaliman and M.\ Zaidenberg proved in  \cite{KZ} 
that every holomorphic mapping $f\colon X\to Y$ from a Stein manifold $X$ to any 
complex manifold $Y$ can be approximated on any compact set $K\subset X$
by holomorphic maps from a neighborhood of $K$ to $Y$ whose $k$-jet extension
is transversal to a given analytic subset of the jet manifold $J^k(X,Y)$. See 
also \cite[Theorem 4.8]{AJM}. This gives the following analogue of 
Theorem \ref{generic1}:

{\em Let $X$ be a Stein manifold, and let $Y$ be a complex manifold. Let $n,p,r,d$ 
be as in Theorem \ref{generic1}. Given a compact set $K\subset X$ and a holomorphic 
map $f\colon X\to Y$, there is a holomorphic map $\wt {f}$ from
an open neighborhood of $K$ in $X$ to $Y$ which approximates $f$ on $K$ as 
close as desired and satisfies
\begin{itemize}
\item[(1)] if $d\leq n$ then $\dim_z \Sigma_{\wt{f},r}= n-d$
for all $z\in K$; 
\item[(2)]  if $d>n$ then ${\rm rank}\; \wt{f} \geq r$ on $K$.
\end{itemize}
}
\end{remark}

% PRELIMINARIES

\section{Preliminaries}
Recall \cite{GH} that a {\em projective algebraic set } (or variety) is a closed subset of a complex
projective space $\P_n$ of the form
\[  
	A= \cap_{j=1}^k \{[z_0\colon \ldots \colon z_n]\in \P_n \colon  p_j(z_0,\ldots,z_n)=0 \}  
\]
where the $p_j$'s are homogeneous holomorphic  polynomials on $\C^{n+1}$. 
Such $A$ is a closed complex analytic subvariety of $\P_n$, and 
every closed complex analytic subvariety of $\P_n$ is of this form by 
Chow's theorem \cite[p.\ 74]{C}.  
Occasionally we shall omit the adjective `projective'. 
The topology on $\P_n$ in which the closed sets are exactly the
projective algebraic sets is called the {\em Zariski topology} on $\P_n$. 
A {\em quasi-projective algebraic set} is a difference $Y\bs Y'$ of two closed 
algebraic subvarieties  $Y,Y'\subset\P_n$. A (quasi) projective {\em algebraic manifold}
is a (quasi-) projective algebraic set without singularities.

Let $U\subset \P_n$ be a quasi-projective algebraic set. 
A function $f\colon U\to \C$ is called a {\em regular function} if 
\[
	f(Z) = f(z_0,\ldots, z_n) =P(z_0,\ldots, z_n)/Q(z_0,\ldots,z_n),
	\quad Z\in U, 
\]
where $P$ and $Q$ are  homogeneous polynomials on $\C^{n+1}$ of the same degree 
and $Q(Z)\neq 0$ for all $Z\in U$. A continuous map $F\colon U\to \P_N$ 
is a {\em regular map} if its components with respect to any
affine chart $\C^N\subset\P_N$ are regular functions on $U\cap F^{-1}(\C^N)$.
If $U'\subset \P_N$ is another quasi-projective algebraic set then
a bijective map $F\colon U\to U'$  is a {\em  biregular isomorphism} 
if both $F$ and $F^{-1}$ are regular maps.

\begin{definition}[{\bf Class $\cA$ manifolds}]
\label{classA}
Let $Y$ be a quasi-projective algebraic manifold.
\begin{itemize}
\item[(i)]
$Y$ is of {\em Class $\cA_0$} if it is covered by finitely many 
Zariski open sets biregularly isomorphic to $\C^m$, $m=\dim Y$.
\item[(ii)] $Y$ is of {\em Class $\cA$} if $Y=\wh{Y}\backslash A$, 
where $\wh{Y}$ is a manifold of class ${\cA_0}$ and $A$ is an
algebraic set in $\wh{Y}$ with complex codimension at least two.
\end{itemize}
\end{definition}

Manifolds of class ${\cA}$ were used by Gromov under the name 
Ell-regular manifolds \cite[\S 3.5]{GR}; our terminology
conforms to the one by Forstneri\v{c} \cite{AJM}.

\begin{example}
Complex affine and projective space, as well as complex Grassmanians,
are manifolds of Class $\cA_0$. Further examples are the 
{\em rational surfaces}, i.e., complex surfaces  birationally equivalent to $\P_2$
\cite[p.\ 244]{BHC}. Apart from $\P_2$ these include the Hirzebruch surfaces 
$\Sigma_n$, $n\in \Z_+$.
\end{example}

% DEFINITION OF TRANSVERSALITY

We recall a few relevant notions regarding transversality of mappings.

\begin{definition}
% [{\bf Transversality}]
Lef $f\colon X\to Y$ be a smooth map of manifolds and let $B$ be a smooth submanifold of $Y$.
We say that $f$ is {\em transverse} to $B$ at a point $x\in X$,
denoted by $f\pitchfork_x B$, if either (i) $f(x)\not \in B$, or (ii) $f(x)\in B$ and 
$T_{f(x)} Y = T_{f(x)} B + df_x (T_x X)$.
If $f\pitchfork_x B$ for all $x\in X$, we say that $f$ is transverse to $B$, 
and denote it by $f\pitchfork B$.
\end{definition}
\rm

For latter application we state a couple of known transverality lemmas. 
The first one provides a lot of transversal maps to choose from 
a transversal family of maps; a proof consists of a reduction
to Sard's theorem and can be found in \cite{AB} or \cite{GP}.

\begin{lemma}[{\bf Transversality Lemma}]
\label{trans}
Let $X,Y,P$ be complex manifolds, let $B$ be a complex submanifold of $Y$, and let 
$\Phi\colon X\times P\to Y$ be a holomorphic map such that 
$\Phi\pitchfork B$. Let $f_t(z) \colon =\Phi(z,t)$. 
Then $\{t\in P \colon  f_t\pitchfork B\}$ is dense in $B$.
\end{lemma}

The following lemma, together with Sard's theorem, implies a jet version of the Transversality Lemma
for holomorphic maps. See \cite[Lemma 4.5]{AJM} or \cite{KZ}; here we supply 
some additional details to the proof.
Recall that $J^k(X,Y)$ denotes the complex manifold of all $k$-jets of 
holomorphic maps $X\to Y$ between a pair of complex manifolds.

\begin{lemma} % [{\bf Jet Transversality Lemma}]
\label{trans1}
Let $X$ be a Stein manifold of dimension $r$, 
embedded as a closed complex submanifold of $\C^n$,
let $Y$ be a complex manifold of dimension $p$, 
and let $F \colon  X\times \C^N\to Y$ be a holomorphic map such that for every 
$x\in X$ the map $F(x,\cdotp)\colon\C^N \to Y$ is a submersion at $0\in\C^N$.
Let $W$ denote the vector space of all holomorphic polynomial maps 
$P\colon \C^n\to \C^N$ of degree $\le k$. For each $P\in W$ set 
$F_P(x)=F(x,P(x))$, $x\in X$. Then the map $H\colon X\times W\to J^k(X,Y)$, defined 
by $H(x,P)=j^k_x(F_P)$, is a submersion in a neighborhood of $X\times \{0\}$
in $X\times \C^N$. 
\end{lemma}

\begin{proof}
We need to prove that $H$ is a submersion at points $(x_0,0)\in X\times W$. Set $F(x_0,0)=y_0$. Choose a neighborhood
$U$ of $x_0$ in $\C^n$ and a neighborhood $V$  of $y_0$ in $Y$ and a neighborhood $E$ of $0$ in 
$W$ such that $\Phi(X\cap U)=\Phi(U)\cap (\C^r\times \{0\}^{n-r})=U'$ for biholomorphic maps 
$\Phi:U\to \Phi(U)\subset \C^n$, $\Psi:V\to \Psi(V)=V'\subset \C^p$ and such
that $F((X\cap U) \times E)\subset V$. We can consider $U'$ an open subset of $\C^r$.
The map $\Psi_1:J^k(X\cap U,V)\to J^k(U',V')$ sending $j_x^k f\in J^k(X\cap U,V)$
to $j^k_{\Phi(x)} (\Psi\circ f\circ \Phi^{-1})$ is well defined and biholomorphic. Let $H'=\Psi_1 \circ H\circ (\Phi|_{U\cap X}\times {\rm id})^{-1}$.
Then $H': U'\times E\to J^k(U',V')$ with
$H'(\Phi(x),P)=j^k_{\Phi(x)} F'(\cdot,P(\cdot))$, where $F'=\Psi\circ F \circ (\Phi|_{X\cap U}\times {\rm id})^{-1}$.
It is enough to show that $H'$ is a submersion at $(\Phi(x),0)$. Therefore we can assume that 
$X=\C^r\subset \C^n$ and also $Y=\C^p$.

Given $x_0\in \C^r$ we denote by $W_{x_0}$ the set of all polynomials $P\in W$ 
such that $P(x_0)=0$. For a fixed $x\in \C^r\subset \C^n$ the map $W_x\to \C^{M(r,N,k)}$,
$P\mapsto \partial^k_x P(x)$, is a submersion.
Here $J^k(\C^r,\C^N)= \C^r\times \C^N\times \C^{M(r,N,k)}$ for some 
$M(r,N,k)\in \N$ and $\partial^k_x P(x)$ denotes all 
partial derivatives of $P$ of order less or equal than $k$ without
the $0$-th derivative. For every multiindex $I=(i_1,\ldots,i_r)$ we can write 
\[
\partial_x^I (F(x,P(x)))=\sum_{j=1}^n \frac{\partial F}{\partial t_j} (x,P(x)) \partial_x^I P_j(x)+R(x).
\]
Here $R$ contains derivatives of $P$ of order lower than $|I|$ and derivatives of $F$.
For a fixed $x=x_0$ we get 
\[
\sum_{j=1}^n \frac{\partial F}{\di t_j} (x_0,0) \di^I_x|_{x=x_0} P_j(x)+R(x_0)= 
\di_t|_{t=0} F(x_0,\cdot) \cdot \di^I_x|_{x=x_0} P(\cdot) +R(x_0),
\]
where $R(x_0)$ depends linearly on the components of $j_{x_0}^{(|I|-1)} P(x)$. We also see that 
$H(x_0,P)$ is a block-wise lower triangular linear map in the base $\{\partial_{x_0}^J P(x), |J|\leq k\}$ of $\C^{M(r,N,k)}$.
Since $\partial_t|_{t=0} F(x_0,\cdot)$ is surjective, the map $P\mapsto j_{x_0}^k(F_P)$ from $W_{x_0}$ to $\C^{N(r,p,k)}$
is a submersion, and hence $H$ is a submersion at $(x_0,0)$.
\end{proof}

\bigskip

\begin{definition}
\label{Ell1}
Let $X$ and $Y$ be complex manifolds. 
Holomorphic maps $X\to Y$ satisfy {\em Condition Ell$_1$} if for every
map $f\in \cH(X,Y)$ there is a holomorphic map $H\colon  X\times \C^N\to Y$ for some $N\in \N$, satisfying
\begin{itemize}
\item[(1)] $H(x,0)=f(x)$ for all $x\in X$, and
\item[(2)]  the map $H(x,\cdotp)\colon \C^N\to Y$ is a submersion at $0\in \C^N$ for every $x\in X$.
\end{itemize}
\end{definition}

The Ell$_1$ condition is useful when combined with the (Jet) Transversality Lemma 
in approximating a given holomorphic map by  a holomorphic map transversal to a given 
submanifold. Condition Ell$_1$ holds for holomorphic maps from any Stein 
manifold to any complex manifold $Y$ which enjoys the CAP property 
\cite[Proposition 4.6 (b)]{AJM}. 
In short the idea is to construct a finite collection of sprays on $Y$ 
using bundles described in Lemma \ref{spray1} and combining them into map 
$H$ from the definition of Ell$_1$. In particular, Ell$_1$ holds for maps
of Stein manifolds to manifolds of Class $\cA$ since these enjoy the CAP property.

The following result from \cite{AJM} follows from Lemma \ref{trans1} and Sard's theorem.

\begin{lemma}
\label{jet_trans1}
{\rm (\cite[Theorem 4.2]{AJM})}
Let $X$ be a Stein manifold and let $Y$ be a complex manifold such that holomorphic maps 
$X\to Y$ satisfy Condition Ell$_1$. Choose a distance function $d$ on $Y$.
Let $Z$ be a closed complex submanifold (or a closed complex subvariety) in $J^k(X,Y)$. 
Given a compact set $K\subset X$, a holomorphic map $f\colon X\to Y$ and 
an $\epsilon>0$, there is a holomorphic map $f_1\colon X\to Y$ 
such that
\begin{itemize}
\item[(1)] $d(f(x),f_1(x)) <\epsilon$ for ever $x\in K$, and
\item[(2)] $j^k f_1\pitchfork Z$ on $K$.
\end{itemize}
\end{lemma}

We will need the following lemma which was also used in the proof of Proposition 2 in \cite{Fo}.

% CALCULATING DIMENSION OF m x n MATRICES OF RANK r

\begin{lemma} 
\label{dimrank}
The set $M^r(n,m)=\{A \in \C^{n\times m} \colon  {\rm rank}\, A = r \}$ is a 
(non-closed) complex submanifold of $\C^{n m}$ of complex codimension 
$(n-r) (m-r)$.
\end{lemma}

\begin{proof} 
Let $A\in M^r(n,m)$. Change bases in $\C^n$ and $\C^m$ such that $A$ takes the form 
$A=\left[\begin{matrix} B & C \\ D & E \end{matrix} \right]$ where $B$ is an invertible 
$r\times r$ matrix. Denote by $U$ the neighborhood of $A$ in $\C^{n m}$ consisting
of all matrices $A'=\left[\begin{matrix} B' & C' \\ D' & E' \end{matrix} \right]$
where $B'$ is invertible. Define a map $\Phi\colon  U\to \C^{(n-r)(m-r)}$   
by $\Phi(A')=E'-D'B'^{-1}C'$. If $B',C',D'$ are fixed, 
this is just a translation, and therefore $\Phi$ is a submersion.

To conclude the proof it now suffices to show $M^r(n,m)\cap U=\Phi^{-1}(0)$.
Let $A'=\left[\begin{matrix} B' & C' \\ D' & E' \end{matrix} \right] \in U$. Note that  
$F=\left[\begin{matrix} I_r & 0 \\ -D' B'^{-1} & I_{n-r} \end{matrix} \right]$ 
is an invertible matrix and hence ${\rm rank}\, A'= {\rm rank}\, (F A')$. 
But $F A' = \left[\begin{matrix} B' & C' \\ 0 & E'-D' B'^{-1} C' \end{matrix} \right]$
which has rank $r$ if and only if $E'-D' B'^{-1} C'=0$, and this is equivalent
to $A' \in \Phi^{-1}(0)$.
\end{proof}

\begin{definition}[{\bf Spray on a manifold}]
A spray on a complex manifold $X$ is a holomorphic map 
$s\colon E\to X$ from total space of a holomorphic vector bundle $p\colon E\to X$
satisfying $s(0_x)=x$ for all $x\in X$. The spray is algebraic if $p\colon E\to X$ 
is an algebraic vector bundle and $s\colon E\to X$ is algebraic map.
\end{definition}

The following lemma is due to Gromov \cite{GR}
(Lemmas 3.5B and 3.5C); see also \cite[Lemma 1.3]{F4}.
Here we supply additional details of the proof.

\begin{lemma}
\label{spray1}
Let $\wh{Y}$ be an $n$-dimensonal manifold of  Class $\cA_0$ and let $U\subset \wh{Y}$ be
a Zariski open subset biregularly isomorphic to $\C^n$ via an isomorphism $\varphi\colon U\to\C^n$.
Let $\Lambda$ be a closed algebraic subset of $\wh{Y}$ of pure dimension $n-1$ such that
$U\cup \Lambda = \wh{Y}$.
Let $s\colon  U\times \C^n \to U$ be a spray defined by 
\[ 
	s(y,t)\colon =\varphi^{-1}(\varphi(y)+t), \quad y\in U,\ t\in\C^n,
\]
and let $L=[\Lambda]^{-1}$ where $[\Lambda]$ is the line bundle defined by the divisor of $\Lambda$.
There are an integer $m\in \N$ and an algebraic spray 
$\tilde{s}\colon  E= (\wh{Y} \times \C^n)\otimes L^{\otimes m}\to \wh{Y}$
such that $\tilde{s}=s$ on $E|_{\wh{Y}\backslash \Lambda}$ and $s(E_y)=\{y\}$ for all $y\in \Lambda$.	
(Here we have  identified $E|_{\wh{Y}\backslash \Lambda}$ with $({\wh{Y}\backslash \Lambda})\times \C^n$ since 
$L|_{\wh{Y}\backslash \Lambda}$ is trivial.)
\end{lemma}

\begin{proof}
We can't just extend $s$ to $\wh{Y}\times\C^n$ because of the singularities on $\Lambda$.
However, we will show that any point $y\in \Lambda$ admits a Zariski neighborhood 
$\C^n \simeq V\subset \wh{Y}$ such that $s$ extends to $E|_V$ for $m\in N$ large enough. 
Since $\wh{Y}=\cup_{j=1}^r U_j$ for Zariski open sets $U_j$ biregularly isomorphic
to $\C^n$, with $U_1=U$, we will get the desired extension by choosing the largest $m$. 

Let $\varphi_j\colon \U_j\to \C^n$, $1\leq j\leq r$ be biregular isomorphisms; 
the collection $\{(U_j,\varphi_j) \colon 1\le j\le r\}$ is then an algebraic atlas 
on $\wh{Y}$. Choose $y_0\in \Lambda\bs U$; without loss or generality we may assume
that $y_0\in U_2$ and $\varphi_2(y_0)=0$. 
Recall that the spray  $s$ is given in the  local chart $U_1\times \C^n$ on $\wh{Y}\times \C^n$
by $s'_1(z,t)=z+t$. In the local chart $U_2\times \C^n$ the same spray is of the form
\[
	s'_2(z,t)= \varphi_{1,2}^{-1}( \varphi_{1,2}(z)+t),
	\quad z\in \varphi_2(U_{1,2})\subset\C^n,\ t\in\C^n,
\]
where $\varphi_{1,2}=\varphi_1\circ \varphi_2^{-1}$ and $U_{1,2}=U_1\cap U_2$.
Clearly $s'_2$ is holomorphic on the set 
\[
	\{(z,t)\in \C^n\times \C^n \colon z\in \varphi_2(U_{1,2}),\ \varphi_{1,2}(z)+t\in \varphi_1(U_{1,2}) \}
\]
and has singularities in the complement.
In particular, $s'_2$ is holomorphic at all points 
$(z,0)$ with $z\in \varphi_2(U_{1,2})$. Since $U_2 \bs U_1\subset\P_n\bs U_1\subset \Lambda$,
we have $\Omega:=\C^n\bs \varphi_2(\Lambda\cap U_2) \subset \varphi_2(U_{1,2})$
and hence $s'_2$ is holomorphic on a neighborhood of $\Omega\times \{0\}$. 

For a fixed $z\in \Omega$ we can write 
$s'_2(z,t)=z+\sum_{|\alpha|=1}^\infty f_{\alpha}(z)t^{\alpha}$, where 
$\alpha$ is a multiindex and $f_\alpha$ are matrices with rational functions as elements.
Note that the transition maps of the bundle $E\to \wh Y$ are 
$\Phi_{ij}: U_{i,j} \times \C^n\to U_{i,j}\times \C^n$ where 
\[
	\Phi_{ij}(y,t)=(y, (b_j(y)/b_i(y))^m t). 
\]
Here $b_j$ is a regular defining function for $\Lambda\cap U_j$ and $b_i$ is a regular 
defining function for $\Lambda\cup U_i$.
The bundle $E|_{U_1}$ is trivial and can be identified with $U_1\times \C^n$. 
Denote by $\tilde{s}'_2$ the map $s$ in the local chart $U_2 \times \C^n$ on $E$.
Then 
\[
	\tilde{s}'_2(z,t)=(s'_1\circ  (\varphi_1\times {\rm id }) \circ \Phi_{12}\circ (\varphi_2 \times {\rm id})^{-1})(z,t)=
	s'_2(z, t \, b_2(\varphi_2^{-1}(z))^m/b_1(\varphi_2^{-1}(z))^m).
\]
For $z\in U_2\cap \varphi_2(U_{1,2})$ this can be written as 
\[
	\tilde{s}'_2(z,t)=z+\sum_{|\alpha|= 1}^\infty f_{\alpha}(z)\cdot t^{\alpha} b_2(\varphi_2^{-1}(z))^ {m |\alpha|}/b_1(\varphi_2^{-1}(z))^ {m |\alpha|}. 
\]
using the above series expansion for $s'_2$.
By the Cauchy formula for the coefficients of a power series for the rational map $s'_2$, holomorphic on a neighborhood of $(z,0)$,
the maximum of the degrees of the poles of $f_\alpha(z)$ is bounded by some integer which is independent of $z$ and $\alpha$. 
Hence there is $m\in \N$
such that $b_2(\varphi_2^{-1}(z))^m f_{\alpha}(z)$ is holomorphic on $\C^n$ and equals zero 
when $z\in \varphi_2(\Lambda\cap U_2)$, $|\alpha|\in \N$.
This shows that for such $m$ the map $\tilde{s}$ is holomorphic on $E$ at points $0\in E_y$, $y\in \Lambda$.
For other points $t\in \C^n$ we can still write $s'_2$ as power series, since the intersection of the singular set of $s'_2$ with 
$\{z\}\times \C^n\equiv \C^n$ is  nowhere dense in $\C^n$. Furthermore, because of the factors $b_2(\varphi_2^{-1}(z))^m$ we 
can extend $\tilde{s}'_2$ to a continuous, locally bounded map on a neighborhood of hypersurface $\varphi_2(\Lambda\cap U_2)\times \C^n$
with $\tilde{s}'_2 (z,t)=z$ for $z\in \varphi_2(\Lambda\cap U_2)$, $t\in \C^n$.
By Riemann extension theorem $\tilde{s}'_2$ extends to a holomorphic map on a neighborhood of $\varphi_2(\Lambda\cap U_2)\times \C^n$.

\end{proof}

% BEGINING OF MAIN PROOFS

\section{Proofs of main theorems}
The following is Lemma 3.4 in \cite[p.\ 156]{F3} with $r=n-2$, $s=2$ and $D\times L$ instead of $L$.

\begin{lemma}
\label{autom2}
Let $K\subset \C^n$ be a product of a closed polydisc and a ball,
and let $\Sigma\subset \C^n\backslash K$ be
an algebraic set  with $\dim \Sigma \leq n-2$. Let $D=\pi(K)$, where $\pi\colon \C^{n-2}\times \C^2\to \C^{n-2}$ 
is a standard projection and $L\subset \C^2$ is a compact 
polydisc such that  $K\subset D\times L$. Given $\epsilon>0$ 
there exists an automorphism $\Psi$ of $\C^n$ of the
form $\Psi(z',z'')=(z',\psi(z',z''))$ $(z'\in\C^{n-2},\ z''\in \C^2)$
such that 
\begin{itemize}
\item[(i)] $|\Psi(z)-z|<\epsilon$ for all $z\in K$, and  
\item[(ii)] $\Psi(D\times L)\subset \C^n\backslash \Sigma$.
\end{itemize}
\end{lemma}

The following lemma is the main ingredient in the proof of Theorem \ref{theorem1}.

% MAIN APPROXIMATION LEMMA

\begin{lemma}
\label{approx1}
Let $Y$ be a manifold of class $\cA$ with $\dim Y=p$. Choose a distance function $d$ on $Y$. Let $K\subset \C^n=\C^{n-2}\times \C^2$, $L\subset \C^{2}$, $D\subset \C^{n-2}$ be as in Lemma \ref{autom2}. Let $r\in \N$ satisfy $(n-r+1) (p-r+1)\ge 2$.
Given a holomorphic map $f\colon K \to Y$ satisfying ${\rm rank} f\geq r$ on $K$ and an $\epsilon>0$, there exists an algebraic map  
$\tilde{f} \colon  D\times L \to Y$ such that 
\begin{itemize}
	\item[(i)] $d(\tilde{f}(z),f(z)) <\epsilon$ for all $z\in K$, and 
	\item[(ii)] ${\rm rank} \tilde{f}\geq r$ on $D\times L$. 
\end{itemize}
\end{lemma}

\begin{proof} 
By Corollary 3.2 in \cite{AJM} we can approximate the map $f\colon K\to Y$ with an algebraic map $\C^n\to Y$.
So we can assume that $f$ is algebraic, defined on whole $\C^n$ and with ${\rm rank}\; f \geq r$ on $K$. Let 
$\Sigma_{f,r} = \{ z\in \C^n \colon  {\rm rank}_z f < r \} $. Then $\Sigma_{f,r} \cap K=\emptyset$
provided the above approximation was good enough. 

If $\dim \Sigma_{f,r} \leq n-2$, Lemma \ref{autom2} furnishes an automorphism 
$\Psi$ of $\C^n$ which approximates the identity on $K$ and satisfies 
$\Psi (D\times  L) \subset \C^n\backslash \Sigma_{f,r}$. The map we are looking for is $\tilde{f}=f\circ \Psi$. 

Now suppose that $\dim \Sigma_{f,r} = n-1$. 
We will reduce this to the previous case $\dim \Sigma_f = n-2$.
This reduction is similiar to the one used in the proof of Proposition 5.4 in \cite{AJM}.
By the definition of a Class $\cA$ manifold we have $Y=\wh{Y}\backslash A$ where $\wh{Y}$ is a manifold of Class $\cA_0$ 
and $A$ an algebraic subset of codimension at least two in $\wh{Y}$.
We will approximate $f$ on $D\times L$ by an algebraic map 
$f_0\colon \C^n\to \wh{Y}$ such that $\dim \Sigma_{f_0,r}\leq n-2$. 
By approximating well enough we will also get 
$f_0(D\times L)\subset \wh{Y}\backslash A$. 
In each step of the approximation a given map $f$ will be replaced by 
a nearby algebraic map $f_1$ such that the corresponding set
$\Sigma_{f_1,r} \subset\C^n$ has less $n-1$ dimensional 
irreducible components than $\Sigma_f$.

Choose a point $z_0\in \Sigma_f$ belonging to exactly one $(n-1)$-dimensional irreducible component 
$\Sigma'$ of $\Sigma_{f,r}$. Set $y_0=f(z_0)$. 
By the definition of a class $\cA$ manifold there is a Zariski open neigborhood 
$U$ of $y_0$ in $\wh{Y}$ which is biregularly isomorphic to $\C^p$, 
where $p=\dim \wh{Y}$. Hence there is a biregular isomorphism 
$\varphi\colon U=\wh{Y}\backslash \Lambda\to \C^n=\P_n\backslash H$ where $H$ is
the plane at infinity in $\P_n$. Since $\varphi$ has poles at $\Lambda$ it can be viewed
as a holomorphic map $\wh{Y}\to \P_n$. Let $\varphi(z_0)\in V\subset \P_n$ where $V\equiv \C^n$.
There is a polynomial $q$ defined on $V$  which vanishes on $\varphi(\Lambda)\cap V$ 
but $q(\varphi(y_0))\neq 0$. The closure of the zero locus of $q$ in $\P_n$ is an algebraic set. 
Denote by $\wh{\Lambda}$ its $\varphi$-preimage. 
Then $\wh{\Lambda}$ is an algebraic set in $\wh{Y}$ of pure dimension $p-1$ such that
$\wh{\Lambda}\cup U=\wh{Y}$ and $y_0\not \in \wh{\Lambda}$.

Let $L$ be a holomorphic line bundle over $\wh{Y}$ defined by the divisor of $\wh{\Lambda}$. 
Using Lemma \ref{spray1} we get a spray  
$s\colon E=(\tilde{Y}\times \C^p)\otimes L^{-m}\to \wh{Y}$ such that $s(x,t)=x+t$ on $U$ (identifying $E|_U$ with
$U\times \C^p$ and $U$ with $\C^p$ via an isomorphism) and $s(x,t)=x$ for all $x\in \wh{\Lambda}$.
Let $\iota \colon  f^\ast E\to E$ be a natural map covering $f$. Here 
$f^\ast E = \{ (z,v) \colon  z\in \C^n, v\in E_{f(z)}\}$ is the pullback of the bundle $p\colon E\to \wh{Y}$.
In local coordinates $\iota$ is just $\iota(z,v)=(f(z),v)$. Since $f$ is algebraic, $f^\ast E$ is an algebraic vector
bundle over $\C^n$. 
By Serre's theorem A \cite{Se} the bundle $f^\ast E$ is generated by finitely 
many algebraic sections, and hence there is a surjective algebraic 
vector bundle map $g\colon  \C^n\times \C^q\to f^\ast E$ for some $q\in \N$. 
We can write $g(z,t)=\sum_{j=1}^q g_j(z) t_j$ where $g_j\colon \C^n\to f^\ast E$ 
are sections. Set $\Phi =s\circ \iota \circ g$,
$\Lambda=f^{-1}(\wh{\Lambda})$, $V=\C^n\backslash \Lambda$. From the above statements it is easy to conclude the
following properties of the algebraic bundle map $H \colon  \C^n\times \C^q \to \wh{Y}$:
\begin{itemize}
\item[--]  $H(z,0)=f(z)$ for $z\in \C^n$,
\item[--] $H(z,t)=z$ for $z\in \Lambda$ and $t\in \C^q$, and 
\item[--] $H(z,\cdot)$ is a submersion on $\C^q$ for every $z\in V$.
\end{itemize}

Let $W$ denote the space of all quadratic polynomial maps $\C^n\to \C^q$.
Set $f_P(z)=H(z,P(z))$ for  $P\in W$.
By Lemma  \ref{dimrank} the set $Z_j=\{(z,y,\alpha)\in J^1(\C^n,\wh{Y}) \colon  {\rm rank}\, \alpha =j \}$ is a
submanifold of $J^1(\C^n,\wh{Y})$ of codimension $(n-j)(p-j)$ and $Z=\cup_{j=0}^{r-1} Z_j$ is a closed subvariety
of codimension $(n-r+1)(p-r+1)\geq 2$. 
By Lemma \ref{trans1} for every $P$ in some open dense subset of $W$ we get
$j^1 f_P\pitchfork Z$ .
If we choose $P$ close to $0$ from subset of polynomials with up to first degree
terms equal zero, we can conclude the following about the map $f_P\colon \C^n\to \wh{Y}$: 
\begin{itemize}
\item[(1)]  $f_P(z)=f(z),\ df_P(z)=d f(z)$ for $z\in \Lambda$,
\item[(2)] $f_P$ approximates $f$ on $D\times L$, and 
\item[(3)] the algebraic set $\Sigma_{f_P,r}=(j^1 f_P)^{-1} (Z)$ has dimension less than 
$n-1$ at every point of $V$.
\end{itemize}

By (3) the only remaining irreducible $n-1$ dimensional components of $f_P$ are those 
lying in $\Lambda$, where they are equal those of $\Sigma_{f,r} \cap \Lambda$. 
Hence the number of $n-1$ dimensional components intersecting with $\Lambda$ has not
increased. At least one component $\Sigma'$ from $\Sigma_{f,r}$ is missing since $z_0\notin \Lambda$. 
Set $\Sigma_1=\Sigma_{f_P,r}$, $f_1=f_P$. The map $f_1 \colon  \C^n\backslash \Sigma_1 \to \wh{Y}$ has ${\rm rank} \geq r$ and 
$\Sigma_1$ has less $(n-1)$-dimensional components than $\Sigma_{f,r}$. 
By repeating this procedure we obtain in finitely many steps the 
desired algebraic map $f_0$ with $\dim \Sigma_{f_0}\leq n-2$.
\end{proof}

\begin{lemma}
\label{approx2}
Let $Y$ be a manifold of class $\cA$ with $\dim Y=p$. Let $K\subset \C^n$ be a product of a 
closed ball and a closed polydisc, and let $Q\subset \C^n$ be a closed polydisc containing $K$. 
Let $r$ be such that $(n-r+1) (p-r+1)\geq 2$. Every holomorphic map $f\colon  K\to Y$ with ${\rm rank}\; f\geq r$
on $K$ can be approximated unikformly on $K$ by algebraic maps 
$\tilde{f} \colon  Q\to Y$ satisfying  ${\rm rank}\, \tilde{f}\geq r$ on $Q$.
\end{lemma}

\begin{proof} 
If $n$ is even, we write $\C^n=\C^2\times \cdots \times \C^2$ 
($n/2$ factors) and let $\pi_j\colon \C^n\to \C^{n-2}$ 
be the projection whose kernel is the $j$-th factor. 
Let $Q=L_1\times \cdots \times L_m$ where $L_j\subset \C^2$ are polydiscs. 
Using Lemma \ref{approx1} we approximate $f=f_0$ on 
$K=K_0$ by $f_1\colon L_1\times \pi_1(K)\to Y$. 
Set $K_1=L_1\times \pi_1(K)$ and approximate 
$f_1$ on $K_1$ by $f_2\colon  L_2\times \pi_2(K_1)\to Y$ 
using  Lemma \ref{approx1} (for purposes of shorter notation the 
coordinates have been permuted). 
By continuing in this fashion we get the desired approximation 
in $m=n/2$ steps. In the case of odd $n$ we use an extra disc.
\end{proof}

{\it Proof of Theorem \ref{theorem1}.} 
By the  definition of a class $\cA$ manifold we have $Y=\wh{Y}\backslash A$,
where $\wh Y$ is a manifold of Class ${\cA}_0$ and $A$ is a closed algebraic subset 
of $\wh{Y}$ of codimension at least two.
Choose a distance function $d$ on $\wh Y$ induced by a complete Riemannian metric.

Suppose that $K$ is a closed polydisc in $\C^n$. 
Choose an exhaustion of $\C^n$ with closed polydiscs $Q_j$, $j\in \Z_+$,
where $Q_0=K\subset Q_1$. 

By Corollary 3.2 in \cite{AJM} we can approximate $f$ 
uniformly on $K$ by an algebraic map $f_0\colon \C^n\to Y$.

By Lemma \ref{jet_trans1} we can approximate $f_0$ uniformly on $Q_0$ by 
a holomorphic map transversal to $A$ on $Q_0$, and hence we can assume 
that $f_0\pitchfork A$ on $Q_0$. Choose a positive number $\delta_0>0$ such that every holomorphic map 
$g\colon \C^n\to \wh{Y}$ with $d(g(z),f_0(z))<\delta_0$ for all $z\in Q_1$ 
satisfies ${\rm rank} g\geq r$ on $Q_0$ and $g\pitchfork A$ on $Q_0$. 
Using Lemma \ref{approx2} and Corollary 3.2 in \cite{AJM}
we get a holomorphic map $f_1\colon \C^n \to Y$ satisying 
${\rm rank} f_1\geq r$ on $Q_1$ and 
$d(f_1(z),f_0(z))<\min (\epsilon/2,\delta_0/2)$ for all $z\in Q_0$.

Proceeding inductively we get a sequence of holomorphic maps 
$f_j\colon \C^n\to Y$ and a decreasing sequence of positive numbers $\delta_j>0$
satisfying the following:
\begin{itemize}
\item[(i)] $d(f_{j+1}(z),f_j(z)) < \min(\epsilon/2^{j+1},\delta_j/2)$ for all $z\in Q_j$ and $j\geq 0$,
\item[(ii)] $\rank f_j\geq r$ on $Q_j$ and $f_j\pitchfork A$ on $Q_j$, and 
\item[(iii)]  every holomorphic map $g\colon \C^n\to \wh{Y}$ with 
$d(g(z),f_j(z)) <\delta_j$ for all $z\in Q_{j+1}$ satisfies $\rank  g\geq r$ on $Q_j$
and $g\pitchfork A$ on $Q_j$. 
\end{itemize}

The sequence of holomorphic maps $f_j$ converges uniformly on the compacts in $\C^n$ to
a holomorphic map $F\colon \C^n\to \wh{Y}$ satisfying $d(F(z),f(z))<\epsilon$ 
for all $z\in Q_0=K$ (a consequence of (i)) and $d(F(z),f_j(z)) <\delta_j$ 
for all $z\in Q_j$ and $j\geq 0$ (because of (ii) and the definition of the numbers
$\delta_j$). 

This implies $F\pitchfork A$ on $\C^n$ and ${\rm rank} F\geq r$ on $\C^n$. To show that 
$F(\C^n)\subset Y=\wh{Y}\backslash A$ suppose $F(z)\in A$ for some $z\in \C^n$. 
The transversality condition $F\pitchfork A$ implies 
${\rm rank}_z F+\dim_z A\geq p=\dim \wh{Y}$. This and 
$F\pitchfork_z A$ implies $g(U)\cap A\neq \emptyset$ for all holomorphic maps
$g\colon U\to \wh{Y}$ close enough to $F$ on some neighborhood $U$ of $z$ in $\C^n$. 
Since $f_j\colon \C^n\to Y$ for all $j\geq 0$, we have a contradiction.
This completes the proof of Theorem \ref{theorem1}.

\medskip

{\it Proof of Theorem \ref{generic1}.} 
{\em Case (1)}:
Let $Z=\cup_{j=0}^{r-1} Z_j$ where $Z_j$ is the submanifold of $J^1(X, Y)$ consisting of all jets 
with ${\rm rank}$ $j$. Now $\Sigma_{f,r}=(j^1 f)^{-1}(Z)$. 
Since the codimension of $Z$ is $(n-r+1)(p-r+1)$ (see Lemma \ref{dimrank}), we will get 
$\dim \Sigma_{f,r}=n-(n-r+1)(p-r+1)$ if $j^1 f\pitchfork Z$. If this holds
then the set $\{ f \in \cH(X,Y)\colon  j^1 f \pitchfork Z \}\subset \Omega$
will satisfy the conclusion of Theorem \ref{generic1}.

Choose an exhaustion of $X$ by compact sets $K_l$, $l\in \N$, and let 
\[
	\Omega_l=\{ f \in \cH(X,Y)\colon  j^1 f \pitchfork Z \; {\rm on}\; K_l \}.
\]
Then $\Omega=\cap_{l\in \N} \Omega_l$. Since $\cH(X,Y)$ is a Baire space, it is enough to show that $\Omega_l$ is
open and dense in $\cH(X,Y)$. Openess follows directly from the definition 
of transversality and the fact that $Z$ is closed (if
we perturb $f$ on neighborhood of $K$ a little, transversality condition will still be satisfied on $K$ by Cauchy
inequality for derivatives). To prove density choose $g\in \cH(X,Y)$ and a compact subset $L\subset X$. Using Lemma 
\ref{jet_trans1} we get a holomorphic map $f$ which approximates $g$ on $L\cup K_l$ and satisfies
the transversality condition in the definition of $\Omega_l$.

{\em Case (2)}:
The assumed inequality implies $\dim X + \dim Z <\dim J^1(X,Y)$. Therefore $j^1 f\pitchfork Z$ implies 
$j^1 f (X)\cap Z=\emptyset$ which is equivalent to ${\rm rank}\, f\geq r$ on $X$. 
The set of such $f$ is open and dense by the proof of case (1).

\end{document}